\documentclass[12pt]{article}
\usepackage{amsfonts, amsmath, graphicx, comment}
\newcommand{\C}{\mathbb C}

\def\Arg{\mathrm{Arg}\, }
\author{Alexandre Eremenko\thanks{Supported by NSF grant DMS-1361836.}$\;$ and
Alexander Fryntov}
\title{Remarks on the Obrechkoff inequality}
\date{\today}
\begin{document}
\maketitle
\begin{abstract}
Let $u$ be the logarithmic potential of a probability measure $\mu$
in the plane
that satisfies 
     $$u(z) \le u(|z|)\quad z\in\C,$$
and $m(t)=\mu\{ z\in\C^*:|\Arg z|\leq t\}$. Then
$$\frac{1}{a}\int_0^a m(t)dt\leq \frac{a}{2\pi},$$
for every $a\in (0,\pi)$. This improves and generalizes
a result of Obrechkoff
on zeros of polynomials with positive coefficients.

2010 AMS Subj. Class. 30C15, 31A15.

Keywords: polynomials, subharmonic functions.
\end{abstract}
\subsection{Introduction}

Distribution of zeros of polynomials with positive coefficients
is an old subject going back to Poincar\'e \cite{P}.
For some recent results we mention \cite{B} and references there.

Obrechkoff \cite{O} proved that for every polynomial
$P$ of degree $d$ with non-negative coefficients, and every
$\alpha\in(0,\pi/2)$, the number of roots
in the sector $\{ z\in\C^*:|\Arg z|\leq \alpha\}$ is at most
$2\alpha d/\pi.$

A general question about distribution of roots
of polynomials with non-negative
coefficients was asked by Subhro Ghosh and Ofer Zeitouni
\cite{Z1} in connection with their research on the large deviation theorems
for zeros of random polynomials \cite{Z2}.

For each polynomial of degree $d$, we consider the
{\em empirical measure} which is a probability measure in the plane
consisting of atoms of charge $m/d$ at every root of multiplicity $m$.
The question of Ghosh and Zeitouni was to describe the closure 
of empirical measures of polynomials with positive coefficients.

Obrechkoff's inequality implies that every measure $\mu$ in this closure
must satisfy 
\begin{equation}\label{obr}
\mu\{ z\in\C^*:|\Arg z|\leq\alpha\}\leq\frac{2\alpha}{\pi},
\end{equation}
for every $\alpha\in(0,\pi/2)$.

A complete description of the closure was given in \cite{BES}.
It is evident, that every polynomial with non-negative coefficients
satisfies
$$|P(z)|\leq P(|z|),$$
and that the empirical measure of $P$ is symmetric with respect to
the real axis.

For every finite measure $\mu$ in the plane we define the potential
\begin{equation}\label{pot}
u_\mu(z)=\int_{|\zeta|\leq 1} \log|z-\zeta|d\mu+\int_{|z|>1}\log\left|
1-z/\zeta\right|d\mu.
\end{equation}

\noindent
{\bf Theorem A.} \cite{BES} {\em A measure $\mu$ belongs to the closure
of empirical measures of polynomials with positive coefficients if and only
if $\mu(\C)\leq 1$, $\mu$ is symmetric with respect to the real axis,
and
\begin{equation}\label{ineq}
u_\mu(z)\leq u_\mu(|z|),\quad z\in\C.
\end{equation}
}
\vspace{.1in}

Theorem A is proved by approximation of arbitrary 
potential satisfying (\ref{ineq}) and $u(z)=u(\overline{z})$
by potentials of the form $\log|P|/\deg P,$
where $P$ is a polynomial
with positive coefficients.

Combining Theorem A with Obrechkoff's inequality
one concludes that
for every finite measure $\mu$, symmetric with respect to the real line,
condition (\ref{ineq}) implies
(\ref{obr}). The proof of theorem A is complicated,
and it is desirable to obtain a direct potential-theoretic proof
of the implication (\ref{ineq})$\rightarrow$(\ref{obr}).
Such a proof will be given in this paper. In fact we will prove a stronger
statement.
\vspace{.1in}

\noindent
{\bf Theorem 1.} {\em Let $\mu$ be a probability measure in the plane,
symmetric with respect to the real line,
whose potential (\ref{pot}) satisfies (\ref{ineq}). Then the function
\begin{equation}\label{1}
m(t)=\mu\{z\in\C^*:0\leq |\Arg z|\leq t\}
\end{equation}
satisfies
\begin{equation}\label{2}
\frac{1}{a}\int_0^a m(t)\, dt\leq\frac{a}{2\pi},\quad 0\leq a\leq\pi.
\end{equation}
}
\vspace{.1in}

For the uniform distribution on the unit circle we have $m(t)=t/\pi$,
and equality holds in (\ref{2}) for all $a$. Obrechkoff's inequality
(\ref{obr}) is an immediate corollary of (\ref{2}): setting $a=2\alpha$,
we obtain
\begin{equation}\label{obb}
m(\alpha)\leq\frac{1}{\alpha}\int_\alpha^{2\alpha}m(t)dt\leq
\frac{2}{a}\int_0^a m(t)dt\leq \frac{a}{\pi}=\frac{2\alpha}{\pi}.
\end{equation}

Next we discuss the possibility of equality in (\ref{obr}).
For the polynomial $P(z)=z^d+1$ with non-negative coefficients
and $\alpha=\pi/d$
we have equality in (\ref{obr}).
Thus (\ref{obr}) is exact for each $\alpha$ of the form $\pi/d,\;
d=2,3,4,5,\ldots$. The second result of this paper is that in fact
(\ref{obr}) is best possible for all $\alpha$.
For each $\alpha\in(0,\pi/2)$
we will find a probability
measure $\mu$ symmetric with respect to the real axis, satisfying
(\ref{ineq}) and such that equality holds in (\ref{obr}). Then it follows
from Theorem A, that the right hand side of (\ref{obr}) cannot be 
replaced by a smaller number if the resulting inequality must
hold for empirical
measures of all polynomials with non-negative coefficients.

However it is not clear whether equality
can hold in (\ref{obr}) 
for an empirical measure of a polynomial when $\alpha$ is rational
but not
of the form $\pi/m$ with integer $m>2$.

\subsection{Proof of Theorem 1}

Without loss of generality we assume that 
the closed support of $\mu$ is bounded and does not contain $0$:
it was shown in \cite{BES} that arbitrary finite measure satisfying
(\ref{ineq}) can be approximated
by a measure with such a support which also satisfies (\ref{ineq}).

Then it sufficient to consider a potential of the form
    $$u(z) : = \int_\C \log|1-z/\zeta|\, d\mu(\zeta)$$
which differs from (\ref{pot}) by an additive constant,
and hence, also satisfies (\ref{ineq}).

For a fixed $\rho\in(0,1)$, consider the function
$$v_\rho(z)=\int_0^\infty u(z/t)t^{\rho-1}dt.$$
This function is subharmonic and homogeneous,
$$v_\rho(\lambda z)=\lambda^\rho v_\rho(z),
\quad\mbox{for every}\quad \lambda>0,$$
therefore it has the form
\begin{equation}\label{v}
v_\rho(re^{i\theta})=r^\rho h_\rho(\theta).
\end{equation}

To relate $h$ with $\mu$, we 
need the integral
$$\int_0^\infty\log\left|1-\frac{z}{t}\right|t^{\rho-1}dt=c_\rho r^\rho\cos\rho(\theta-\pi),\quad z=re^{i\theta},\; 0\leq\theta\leq 2\pi,$$
where $c_\rho=\pi/(\rho\sin\pi\rho)$. Let us define $\phi_\rho$
as the $2\pi$-periodic extension of
$\cos\rho(\theta-\pi),\;0\leq\theta\leq2\pi.$
Then we have
$$v_\rho(re^{i\theta})=\int_0^\infty\int_\C
\log\left|1-\frac{re^{i\theta}}{t\zeta}\right|
d\mu(\zeta)t^{\rho-1}dt
=c_\rho r^\rho\int_\C\phi_\rho(\theta-\arg\zeta)\frac{d\mu(\zeta)}{|\zeta|^\rho}$$
Comparing this with (\ref{v}) we obtain
\begin{equation}\label{h}
h_\rho(\theta)=\int_0^{2\pi}\phi_\rho(\theta-t)d\nu_\rho(t),
\end{equation}
and 
\begin{equation}\label{nu}
\nu_\rho(E)=c_\rho\int_{\zeta/|\zeta|\in E}\frac{d\mu(\zeta)}{|\zeta|^\rho},
\end{equation}
for every Borel set $E$ on $[0,2\pi)$.
When $\rho\to 0$, $\nu_\rho/c_\rho\to\nu_0$, where $\nu_0$ is proportional to
the
radial projection of the measure $\mu$, so $m(t)=\nu_0[-t,t]).$

Inequality (\ref{ineq}) and symmetry $u(z)=u(\overline{z})$ imply
\begin{equation}
2h_\rho(0)-h_\rho(a)-h_\rho(-a)\geq 0,\quad a\in[0,\pi].
\end{equation}
Using the expression (\ref{h}) we conclude that 
$$\int_0^{2\pi}J_\rho(t)d\nu_\rho(t)\geq 0,$$
where
$$J_\rho(t)=2\phi_\rho(t)-\phi_\rho(t-a)-\phi_\rho(t+a).$$
Now we divide by $\rho^2$ and pass to the limit 
$\rho\to 0$, using $\cos t\sim 1-t^2/2$.
A simple direct computation shows that $J_\rho/\rho^2\to J$,
where
$$J(t)=\left\{\begin{array}{ll}4\pi |t|-4\pi a+2a^2,&|t|\leq a,\\
&\\
2a^2,&a<|t|\leq\pi.\end{array}\right.$$

We conclude that
$$\int_{-\pi}^\pi J(t)\, d\nu_0(t)\geq 0,\quad\mbox{and thus}\quad
\int_{-\pi}^\pi J(t)\, d\mu(t)\geq 0.$$
Integrating the last integral by parts, we obtain
$$4\pi\int_0^a m(t)\, dt \leq J(\pi)m(\pi)=2a^2,$$
which is equivalent to (\ref{2}).

\subsection{Example}

In this section, for any given $\alpha\in(0,\pi/2)$ we construct a probability
measure $\mu$ symmetric with respect to the real line,
and satisfying (\ref{ineq}), such that Obrechkoff's inequality (\ref{obr}) holds with equality.

Inequalities (\ref{obb}) suggest that the sectors
$|\Arg z|<\alpha$ and $|\Arg z|\in(\alpha,2\alpha)$ must be free of
the measure.

Potential
    
$$u(z) : = \log |z^2 + 1|$$
satisfies (\ref{ineq}), and its total Riesz' measure equals $2$\,.
Take $\alpha\in (0,\pi/2)$ and  define the subharmonic function
$$u_\alpha(z) := 
    \left\{
    \begin{array}{ll}
    u(z^{\pi/(2\alpha)}), & |\Arg(z)|< 2\alpha,\\
&\\
    u(|z|^{\pi/(2\alpha)}), & \text{otherwise.}
    \end{array}
    \right.
$$
It is clear that $u_\alpha$ satisfies (\ref{ineq}).
Let $\lambda_\alpha$ be the Riesz' measure of $u_\alpha$.
One should notice that $\lambda_\alpha$ is supported on the set
$$\{z:\,|\Arg(z)|\ge 2\alpha\}\cup\{e^{i\alpha}\}
    \cup\{e^{i\alpha}\}\,.$$
Notice that $\lambda_\alpha\{e^{\pm i\alpha}\} = 1$, and
$\lambda_\alpha$ is absolutely continuous on 
\hbox{$\{z:\,|\Arg(z)|\ge 2\alpha\}$} with respect to
the plane Lebesgue measure
and its density is
$$     \rho_\alpha = \frac{1}{2\pi}\Delta u_\alpha.$$
Since $u_\alpha(e^{i\theta})$ does not depend on $\theta$
for $|\theta|\in(2\alpha,\pi)$,
we compute the Laplacian $\Delta u_\alpha$ in polar coordinates 
($z= re^{i\phi}$) as follows:
$$    \rho_\alpha(r^{i\phi}) =
    \frac{1}{2\pi r}\,\frac{\partial}{\partial r}
    \left( r\frac{\partial u_\alpha}{\partial r}\right) = 
    \frac{1}{2\pi r}\, \frac{d}{dr}\,\left(
    r\,\frac{d}{dr}\,\log(1+r^{\pi/\alpha})
    \right)
    =
    \frac{1}{2\alpha r}\,\frac{d}{dr}
    \left( \frac{ r^{\pi/\alpha}}{1+r^{\pi/\alpha}}\right). 
    $$
Thus, 
    $$
    \lambda_\alpha\{z:|\Arg(z)|\ge 2\alpha\}) =
    \frac{(2\pi - 4\alpha)}{2\alpha}\, \int_0^\infty
    r\rho_\alpha(r)\, dr = 
    \frac{\pi -2\alpha}{\alpha}\,, 
    $$   
and
    $$
    \lambda_\alpha\{\C\} = 2 + \frac{\pi - 2\alpha}{\alpha} = 
    \frac{\pi}{\alpha}\,.
    $$
Then we define normalized measure
$\mu_\alpha := \lambda_\alpha /\lambda_\alpha(\C)$, 
and
    $$
    \mu_\alpha\{e^{\pm i\alpha}\} = \frac{\alpha}{\pi}.
    $$
So the measure $\mu_\alpha$ satisfies the equation
    $$
    \mu_\alpha\{|\Arg(z)|\le  \alpha \} = \frac{2\alpha}{\pi}.
    $$

\end{document}